\documentclass[11pt]{article}

\usepackage{amsthm}
\usepackage[margin=30mm]{geometry}
\usepackage{amsmath}%
\usepackage{amsfonts}%
\usepackage{amssymb}%
\usepackage{amscd}
\usepackage{graphicx}
\usepackage{booktabs}
\usepackage{color}
\usepackage{xcolor}
\usepackage{bbm}
\usepackage[all,cmtip]{xy}
\usepackage{rotating}
\usepackage{currfile}
\usepackage{enumerate}
\usepackage{subcaption}
\usepackage{enumitem}
\usepackage{wrapfig}
\usepackage{eurosym}
\usepackage{epigraph}
\usepackage{courier}
\usepackage[utf8]{inputenc}
\usepackage[T1]{fontenc}
\usepackage[varqu,varl]{zi4}
\definecolor{grau}{rgb}{0.2,0.2,0.2}
\usepackage[colorlinks=true, urlbordercolor=grau, linkcolor=grau, citecolor=grau, urlcolor=grau]{hyperref}

\makeatletter
\ifnum\inputlineno=\m@ne
\let\showlineno\@empty
\else
\def\showlineno{\the\inputlineno}
\fi
\makeatother
\graphicspath{{img/}}

\numberwithin{equation}{section}

\newtheorem{theorem}[equation]{Theorem}

\newtheorem{lemma}[equation]{Lemma}

\newtheorem{conjecture}[equation]{Conjecture}

\newtheorem{question}[equation]{Question}

\theoremstyle{definition}

\newtheorem{remarks}[equation]{Remarks}

\newcommand{\QQ}{\mathbbm{Q}} 
\newcommand{\NN}{\mathbbm{N}} 
\newcommand{\ZZ}{\mathbbm{Z}} 

\newcommand{\textdef}[1]{\textnormal{\textit{#1}}}

\newcommand{\eps}{\varepsilon}

\newcommand{\wt}[1]{\widetilde{#1}}

\newcommand{\rad}{\textnormal{rad}}

\newcommand{\Qbar}{\overline\QQ}

\newcommand{\OS}{\mathcal{O}_S}
\newcommand{\divides}{\,|\,}

\renewcommand{\phi}{\varphi}

\newcommand{\ord}{\textnormal{ord}}
\newcommand{\noopsort}[1]{}

\newcommand{\gt}{>}
\newcommand{\amp}{&}

\DeclareFontFamily{U}{wncy}{}
\DeclareFontShape{U}{wncy}{m}{n}{<->wncyr10}{}
\DeclareSymbolFont{mcy}{U}{wncy}{m}{n}
\DeclareMathSymbol{\Sh}{\mathord}{mcy}{"58} 
\newcommand{\Sha}{\Sh}
\newcommand{\email}[1]{#1} 

\begin{document}

\title{Elliptic curves with good reduction\\ outside of the first six primes}

\author{%
Alex J. Best%
\thanks{%
    Boston University, \email{alex.j.best@gmail.com}. 
Supported by Simons Foundation grant \#{}550023 and a Hariri Institute Graduate Student Fellowship.%
} 
 \and Benjamin Matschke%
\setcounter{footnote}{6}%
\thanks{%
Boston University, \email{matschke@bu.edu}. 
Partially supported by Excellence Initiative of Universit\'e de Bordeaux (IdEx project DiGeMANT), and by Simons Foundation grant \#{}550023.%
} \\ ~
\setcounter{footnote}{-1}%
}
\date{July 20, 2020}




\maketitle


\begin{abstract}
We present a database of rational elliptic curves, up to $\QQ$-isomorphism, with good reduction outside $\{2,3,5,7,11,13\}$. We provide a heuristic involving the $abc$ and BSD conjectures that the database is likely to be the complete set of such curves.
Moreover, proving completeness likely needs only more computation time to conclude.
We present data on the distribution of various quantities associated to curves in the set.
We also discuss the connection to $S$-unit equations and the existence of rational elliptic curves with maximal conductor.
\end{abstract}

%
%
%

\section{Introduction}


Databases or tables of all elliptic curves subject to various constraints have been published since the 1970s, including in the well known Antwerp IV conference proceedings~\cite{BirchKuyk}. 
Such tables are useful both in identifying a given curve appearing in nature, or for proving a curve with certain properties does not exist. Tables can also be used to answer distributional questions about properties of elliptic curves when ordered in different ways.

The most well known such tables are those of elliptic curves over \(\QQ\) with bounded conductor due to Cremona~\cite{Cremona97book,CremonaData}, which now form part of the LMFDB~\cite{lmfdb}.

One may instead however construct tables of elliptic curves with bad reduction only at primes in a specified set of rational primes $S$. These are exactly the primes dividing the conductor.
Organising curves by their primes of bad reduction can be quite useful in practise; it is often possible to prove a particular curve has good reduction outside certain places, and then conclude that the curve is contained in such a table for some $S$.

In particular many classical diophantine equations can be phrased in terms of the existence of elliptic curves with specified places of bad reduction, see Sections~\ref{secSUE_connection} and~\ref{secReductionsOfDiphantineEquations}.



\begin{samepage}
In this paper we compute and study what is conjecturally the complete set of isomorphism classes of elliptic curves over \(\QQ\) with good reduction away from the first six primes \(\{2,3,5,7,11,13\}\).
This set and the code and auxiliary data used to compute it (including Mordell--Weil bases for almost 100,000 Mordell curves) are available at\begin{center}
\href{https://github.com/elliptic-curve-data/ec-data-S6}{\textit{https://github.com/elliptic-curve-data/ec-data-S6}}.
\end{center}
\end{samepage}
Many of the curves in this set have quite large conductor, but nevertheless by virtue of having bad reduction at only a few small primes can be simpler arithmetically than other curves with smaller conductor.

\paragraph{History.}
We now give a non-exhaustive overview of previous work computing databases of elliptic curves over~$\QQ$.

In the late 1980's Brumer and McGuinness~\cite{brumermcguinness} computed rational elliptic curves of prime discriminant bounded by $|\Delta|\leq 10^8$. 
Stein and Watkins~\cite{SteinWatkins} then extended this database to include almost all curves up to $|\Delta|\leq 10^{12}$ with either conductor $N\leq 10^8$ or prime conductor less than~$10^{10}$.

To compute the set of elliptic curves with with bounded conductor, Tingley~\cite{Tingley75thesis} used modular symbols to find all elliptic curves with $N\leq 200$.
This was greatly extended and improved by Cremona~\cite{Cremona97book,CremonaData}, who has currently computed all of these curves up to $N\leq 500000$.
Initially this approach was only known to compute modular elliptic curves, and it was only when modularity was proved that it was confirmed~\cite{BCDT01modularityOverQ} that over~$\QQ$ being modular is not a restriction.


A third natural basis on which to construct a database of elliptic curves, is by restricting the set of places of bad reduction, i.e.\ the primes that divide~$N$ (or equivalently, primes that divide the minimal discriminant).
For any finite set of rational primes~$S$, let $M(S)$ denote the finite set of elliptic curves over $\QQ$ with good reduction outside of~$S$, up to $\QQ$-isomorphism, and let
\[
N_S:=\prod_{p\in S} p\text.
\]
We may then hope to compute the set $M(S)$ for various sets $S$.

The set $M(\{2,3\})$ was computed by Coghlan~\cite{Coghlan67ellipticCurves23} and Stephens~\cite{Stephens65thesis}, and Coghlan's data was republished as Table 4 in~\cite{BirchKuyk}. 
Agrawal, Coates, Hunt and van der Poorten~\cite{AgrawalCoatsHuntVDPoorten80} computed $M(\{11\})$ via a reduction to Thue--Mahler equations.
Cremona and Lingham~\cite{CremonaLingham07ellipticCurves} computed $M(\{2,p\})$ for $p\leq 23$ via a reduction to the computation of $S$-integral points on Mordell curves.
Koutsianas~\cite{Koutsianas19ellipticCurvesOverNFs} used a reduction to $S$-unit equations over number fields to compute $M(\{2,3,23\})$, as well as curves $E\in M(S)$ for various other $S$ satisfying certain restrictions on the $2$-division field of~$E$.
Von K\"anel and the second author~\cite{vKMa14sUnitAndMordellEquationsUsingShimuraTaniyama} computed $M(\{2,3,5,7,11\})$ as well as all $M(S)$ with $N_S \leq 1000$ using an elliptic logarithm sieve to compute $S$-integral points on elliptic curves.
Bennett and Rechnitzer~\cite{BennettRechnitzer} and
Bennett, Gherga and Rechnitzer~\cite{BennettGhergaRechnitzer19ellipticCurvesOverQ} computed $M(\{p\})$ for all $p\leq 50000$ using a refinement of the reduction to Thue--Mahler equations and Thue equations. The latter paper also recomputes $M(\{2,3,5,7,11\})$ using this approach.
Moreover, using a heuristic they computed all curves in $M(\{p\})$ for $p\leq 10^{10}$, without guaranteeing completeness.

Finally we mention that there are various extensions to the above methods to compute elliptic curves over number fields with good reduction outside of a given set of places. In particular the aforementioned approaches of Cremona and Lingham~\cite{CremonaLingham07ellipticCurves} and of Koutsianas~\cite{Koutsianas19ellipticCurvesOverNFs} generalise to the number field setting.

\medskip

\paragraph{Outline.}
The aim of this paper is to compute the set $M(\{2,3,5,7,11,13\})$.
We have computed a subset of this is heuristically the full set, but is not proved to be complete by our method at present.\footnote{However work in progress by the second author gives the same set of curves using a different method.}
In Sections~\ref{secSummary} and~\ref{secDistribution} we give a summary of our data and discuss some statistics of the data.
We compare our data to Cremona's database in Section~\ref{secComparisionToCremonaDB}.

Our computation relies on a reduction to solving Mordell equations in $S$-integers, this is discussed in Section~\ref{secComputationMethod}.
The main computational bottleneck is to compute the Mordell--Weil bases of a large set of Mordell curves, this is elaborated upon in Section~\ref{secMW}.

In Sections~\ref{secCompletenessHeuristic} and~\ref{secSHallandABC} we discuss a heuristic that our database should be complete, and the possibility of proving completeness via additional computation.
In Section~\ref{secMaximalConductor} we show some results suggested by the data regarding the question for which sets $S$ there are elliptic curves with good reduction outside $S$ of maximal possible conductor.
In Section~\ref{secApplications} we discuss connections and applications to solving other classical diophantine equations including $S$-unit, Thue--Mahler and Ramanujan--Nagell equations.

\paragraph{Acknowledgement.}
It is our pleasure to thank Edgar Costa for various useful comments and for computing the analytic ranks of all curves in our database, as well as the leading coefficients and root numbers of the associated $L$-series.
They are available from the same GitHub repository.

\subsection{Summary of the database}
\label{secSummary}
Let $S(n)$ denote the set of the first $n$ rational primes.
According to our computation, the set $M(S(6))$ contains $4576128$ curves in total; see Table~\ref{table:counts}.
Here, $j(M(S(n)))$ is the set of distinct $j$-invariants of curves in $M(S(n))$, the cardinality of this set is therefore the number of $\Qbar$-isomorphism classes of curves in $M(S(n))$.

\begin{table}[!ht]
\centering
\begin{tabular}{lllr}
\toprule
$n$ & $\# M(S(n))$ & $\#j(M(S(n)))$\tabularnewline
\midrule
$0$ & $0$         & $0$       & Tate (cf.\ Ogg~\cite{Ogg66_2powerConductor}) \tabularnewline
$1$ & $24$        & $5$       & \cite{Coghlan67ellipticCurves23,Stephens65thesis,Ogg66_2powerConductor} \tabularnewline
$2$ & $752$       & $83$      & \cite{Coghlan67ellipticCurves23,Stephens65thesis} \tabularnewline
$3$ & $7600$      & $442$     & \cite{vKMa14sUnitAndMordellEquationsUsingShimuraTaniyama} \tabularnewline
$4$ & $71520$     & $2140$    & \cite{vKMa14sUnitAndMordellEquationsUsingShimuraTaniyama} \tabularnewline
$5$ & $592192$    & $8980$    & \cite{vKMa14sUnitAndMordellEquationsUsingShimuraTaniyama, BennettGhergaRechnitzer19ellipticCurvesOverQ} \tabularnewline
$6$ & $4576128^*$ & $34960^*$ & this paper \tabularnewline
\midrule
\end{tabular}
\caption{Numbers of elliptic curves with good reduction outside $S(n)$ up to $\QQ$-isomorphism and up to $\Qbar$ isomorphism.
The asterisk refers to the possible incompleteness of this paper's table.
The case $n=0$ is the classical result that there is no elliptic curve over $\QQ$ with everywhere good reduction.
}
\label{table:counts}
\end{table}%

When $n \ge 2$ we can obtain all of $M(S(n))$ by taking a representative of each $\Qbar$-isomorphism class of curves in $M(S(n))$ and twisting this representative by all integers divisible only by primes in $S(n)$.
For $j \ne 0,1728$ we only have quadratic twists, when $j = 1728$ we have quartic twists, and for $j= 0$ sextic twists (our assumption that $n\geq 2$ implies that $0,1728 \in j(M(S(n)))$), giving the equation
\[
\#M(S(n)) = 2^{n+1}(\#j(M(S(n))) - 2) + 2\cdot 4^{n} + 2\cdot 6^n\text.
\]
This holds in all cases above, and provides a quick check that nothing that obviously should have be in the database has been missed.


Each curve in $M(S(6))$  has conductor $N\divides 2^83^55^27^211^213^2$, which gives, a priori, $4374$ possibilities for~$N$.
It turns out that exactly $4344$ of them are indeed attained by curves in our set.
The $30$ exceptions for which there is no curve with that conductor are
\[
\begin{split}
\{& 1,  2,   3,   4,   5,   6,   7,   8,   9,   10,  12,  13,  16,  18,  22,  25,  28,  60,  81, 165, \\ &  169, 351, 945, 1280,1820,2673,2816,9984,13365,362880\}.
\end{split}
\]
These exceptions factor as follows
\[
\begin{split}
    \{& 1, 2, 3, 2^{2}, 5, 2 \cdot 3, 7, 2^{3}, 3^{2}, 2 \cdot 5, 2^{2} \cdot 3, 13, 2^{4}, 2 \cdot 3^{2}, 2 \cdot 11, 5^{2}, 2^{2} \cdot 7, 2^{2} \cdot 3 \cdot 5, 3^{4}, 3 \cdot 5 \cdot 11, \\ & 13^{2}, 3^{3} \cdot 13, 3^{3} \cdot 5 \cdot 7, 2^{8} \cdot 5, 2^{2} \cdot 5 \cdot 7 \cdot 13, 3^{5} \cdot 11, 2^{8} \cdot 11, 2^{8} \cdot 3 \cdot 13, 3^{5} \cdot 5 \cdot 11, 2^{7} \cdot 3^{4} \cdot 5 \cdot 7 \}.
\end{split}
\]
These (non-)conductors are all within the range of Cremona's database, and we can therefore check that there are indeed no elliptic curves with any of these numbers as their conductor.
We note that the largest conductor for which no elliptic curve of that conductor exists is less than the square root of the largest possible conductor of a curve in~$M(S(6))$.


\medskip

Next we consider isogeny classes in~$M(S(6))$.
This is also a natural partition of curves in the database as $M(S(n))$ is closed under taking isogenies (any two isogeneous curves have the same conductor).
Our data contains $3688192$ disjoint isogeny classes in total: $2966912$ classes of cardinality~$1$, $646784$ of cardinality~$2$, $4608$ of cardinality~$3$, $60928$ of cardinality~$4$, $6784$ of cardinality~$6$, $2176$ of cardinality~$8$, and no others.
%
%
An example of a curve in $M(S(6))$ with isogeny class of cardinality~$8$ is 
\[
y^2 = x^3 + 827614112325\,x + 276113445805174250.
\]

Edgar Costa has computed the analytic ranks of all curves in our table, as well as the leading coefficients and roots numbers of the associated $L$-series.
His computations use interval arithmetic and hence the leading coefficients are given with exact error bounds.
The standard problem that remains is that it is impossible to verify numerically that the lower derivatives vanish exactly, and thus the computed analytic rank is actually only an upper bound once the rank is large enough.  
According to his computations, there are $1884428$ curves of analytic rank~$0$ in our data, $2267261$ of analytic rank~$1$, $406309$ of analytic rank~$2$, $18003$ curves of analytic rank~$3$, and the remaining $127$ curves are of analytic rank~$4$.
We can compare this to the number of rational elliptic curves with conductor bound $N\leq 500000$ with each rank using Cremona's database:
For these curves, Cremona computed analytic and algebraic ranks (and checked that they coincide), and found that there are $1632686$ curves of rank~$0$, $2124004$ of rank~$1$, $461670$ of rank~$2$, $11243$ of rank~$3$, and $1$ of rank~$4$.
In both tables, we observe a similar larger number of rank~$1$ curves than rank~$0$ curves.
An intriguing difference is the larger number of rank~$4$ curves in our data, compared to a similar total number of curves when ordered by conductor.

\subsection{Distribution of quantities}
\label{secDistribution}
In this section we study the distribution of various arithmetical quantities associated to curves in  our dataset.
As these curves have bad reduction at only the first six primes, they are quite structured and it is interesting to compare answers to distributional questions to when curves are ordered with respect to conductor or discriminant.


\begin{figure}[thb]
	\centering
	\begin{minipage}[b]{0.47\textwidth}
		\centering
		\includegraphics[width=\textwidth]{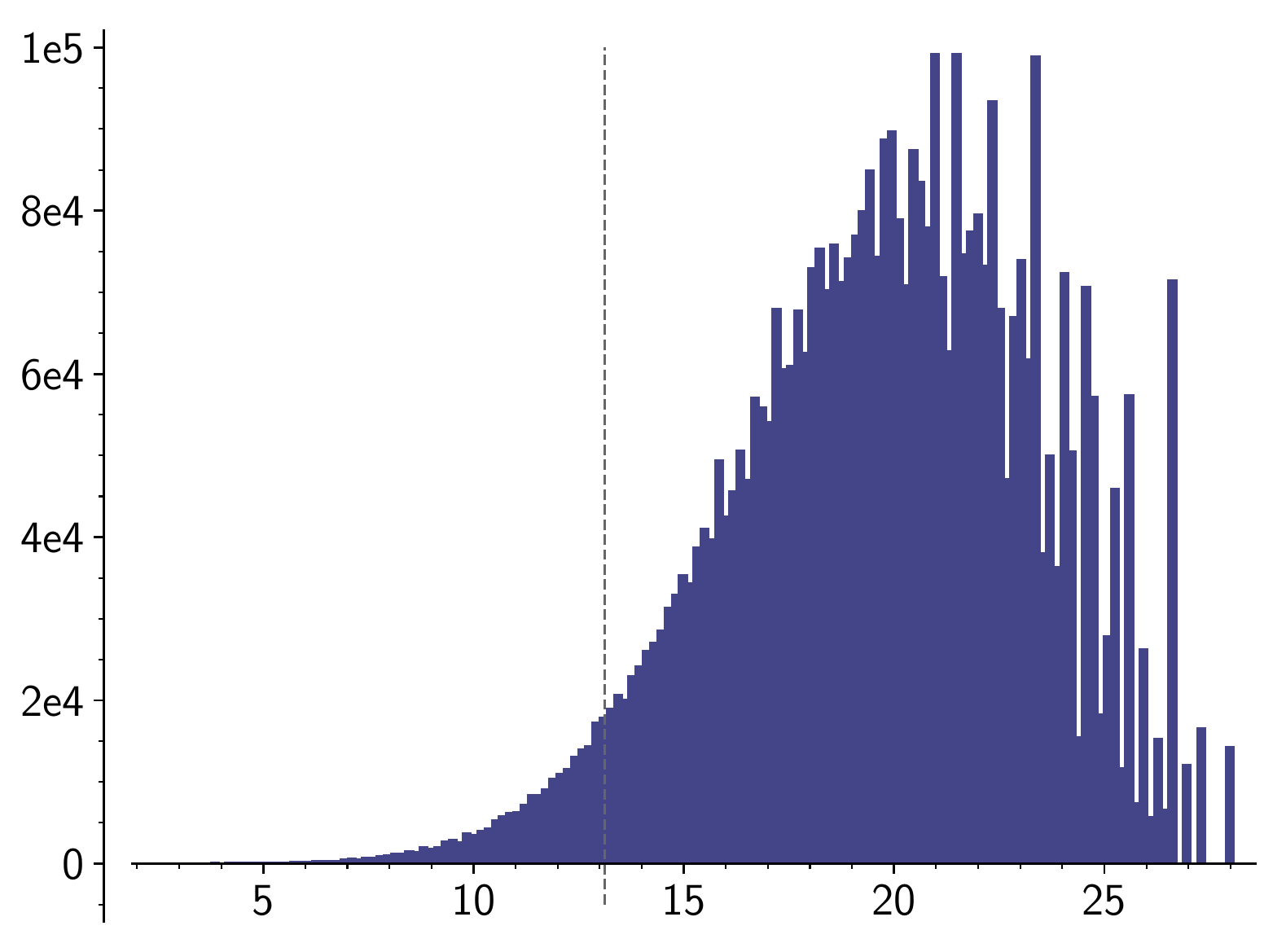}
		\subcaption{Elliptic curves in our data set.}
	\end{minipage}
	\hfill
	\begin{minipage}[b]{0.47\textwidth}
		\centering
		\includegraphics[width=\textwidth]{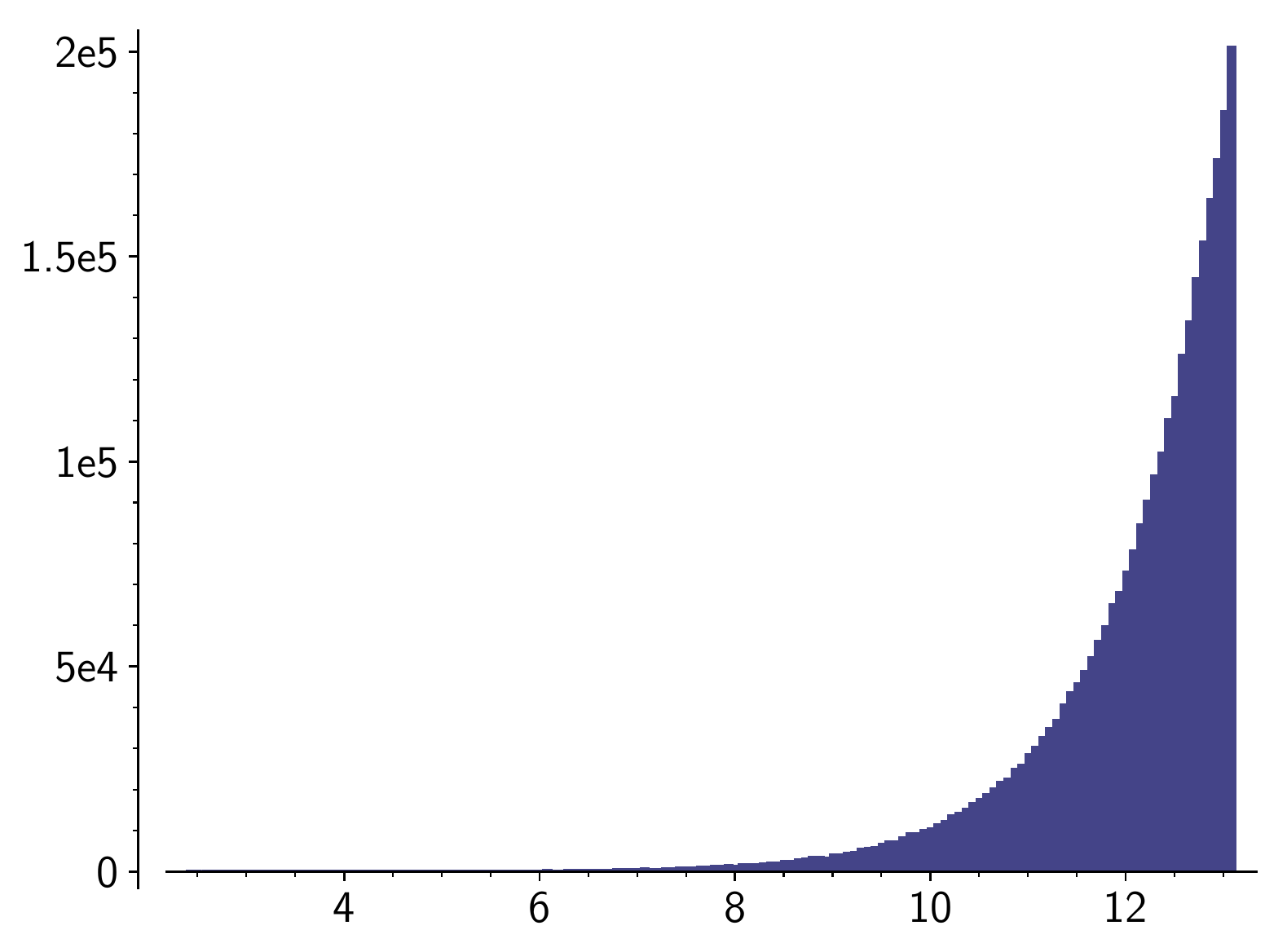}
		\subcaption{Elliptic curves with $N \leq 500000$.}
	\end{minipage}
	\caption{Histograms of logarithms of conductors:
		(a) shows the curves we computed within~$M(S(6))$.
		For a comparison, (b) shows all rational elliptic curves with $N\leq 500000$ according to Cremona's database.
		The bar at $\log(500000)\approx 13.1$ signifies the end of the overlap of both tables.}
	\label{figLogConductors}
\end{figure}

One fundamental quantity is the conductor. 
We plot the distribution of the logarithm of the conductor for the curves in our database as a histogram in Figure~\ref{figLogConductors}(a).
We take the logarithm of $N$ due to the multiplicative nature of the conductor.
Indeed, if the conductor exponents $f_p$ in $N=\prod_{p\in S}p^{f_p}$ were uniformly and independently distributed (which they are not), then in Figure~\ref{figLogConductors}(a) we would see an approximately normal distribution with mean~$14.037$ and standard deviation $4.382$.
The observed distribution of $\log(N)$ is comparatively lopsided: It appears denser in the larger conductor range.
This could be explained by the fact that one can turn good into additive reduction at $p\geq 3$ via twisting by~$p$ (as the reduction of $E$ at $p$ will have Kodaira symbol~$\mathrm{I}_0^*$ by Tate's algorithm), without leaving $M(S(6))$.


\begin{figure}[thb]
	\centering
	\begin{minipage}[b]{0.47\textwidth}
		\centering
		\includegraphics[width=\textwidth]{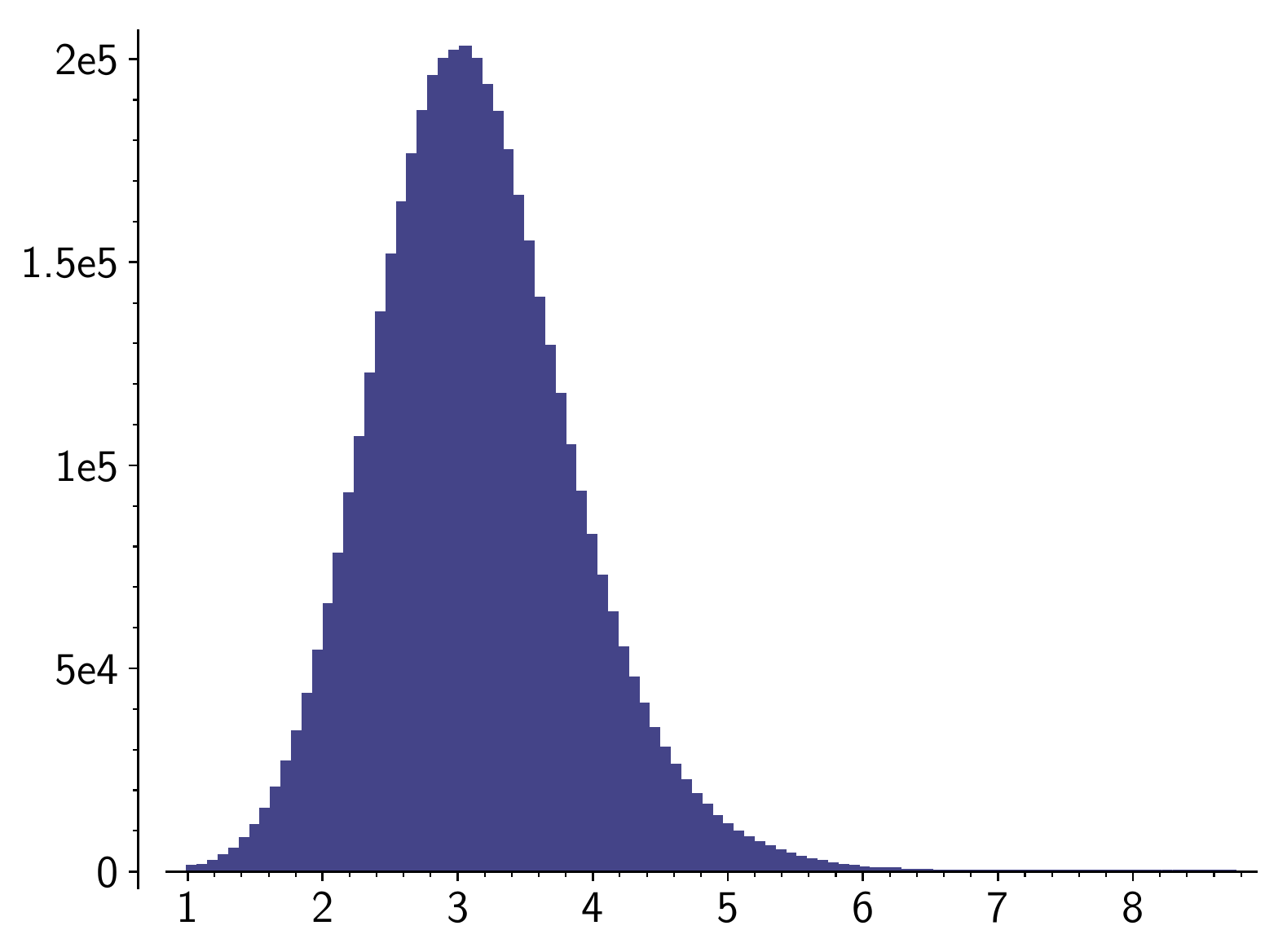}
		\subcaption{Elliptic curves in $M(S(6))$.}
	\end{minipage}
	\hfill
	\begin{minipage}[b]{0.47\textwidth}
		\centering
		\includegraphics[width=\textwidth]{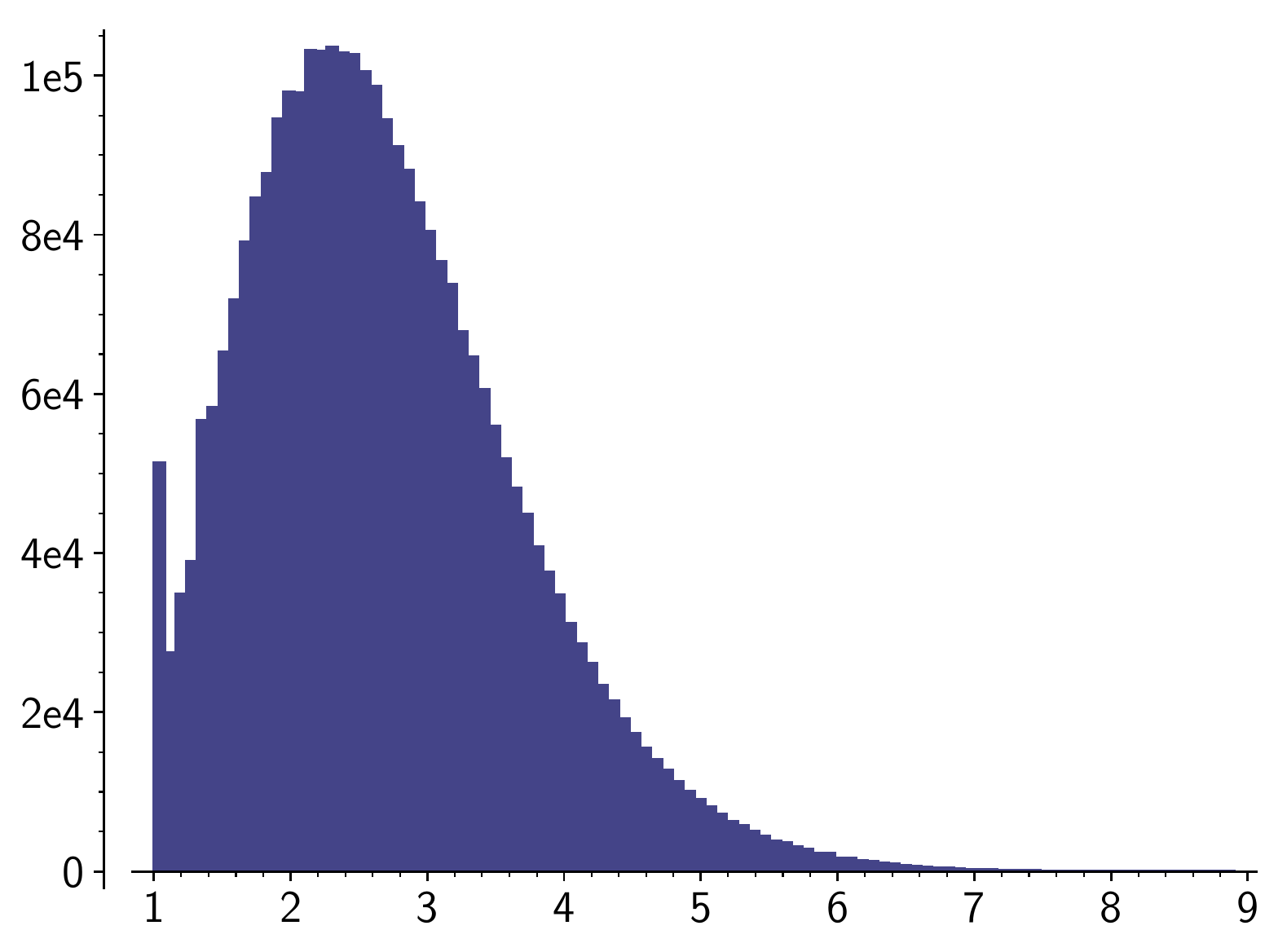}
		\subcaption{Elliptic curves with $N \leq 500000$.}
	\end{minipage}
	\caption{Histograms of Szpiro ratios $\sigma=\log(\Delta_E)/\log(N)$.
		(a) shows the curves we computed within~$M(S(6))$.
		For a comparison, (b) shows all rational elliptic curves with $N\leq 500000$ according to Cremona's database.
	We observe three differences: 
	the larger maximal value for~$\sigma$ in~(b) (namely $8.903700$), 
	the larger mean for~$\sigma$ in~(a), and that (b) contains a significant number of curves with~$\sigma=1$ (namely $602$). 	
	}
	\label{figSzpiroRatios}
\end{figure}

The \emph{Szpiro ratio} of an elliptic curve over $\QQ$ is defined to be the ratio
\[\sigma= \frac{\log|\Delta_{E}|}{\log N}\]
of the logarithms of the minimal discriminant of the curve and its conductor.
Figure~\ref{figSzpiroRatios}(a) sketches the distribution of Szpiro ratios for the curves in our database. 
Szpiro's conjecture states that $| \Delta_E| =  O_\varepsilon  (N^{6 + \varepsilon})$, or equivalently that for any $\delta \gt 0 $ there are only finitely many elliptic curves over $\QQ$ with $ \sigma \gt 6  + \delta$.

Indeed the largest Szpiro ratios occurring among all curves in our dataset are approximately
\[8.757316,\, 8.371586,\, 8.11481\text{ and } 8.034917\ldots\text,\]
the curves for which these ratios occur are all in the LMFDB and have labels \href{https://lmfdb.org/EllipticCurve/Q/858.k2}{\textit{858.k2}}, \href{https://lmfdb.org/EllipticCurve/Q/2574.j2}{\textit{2574.j2}}, \href{https://lmfdb.org/EllipticCurve/Q/910.e1}{\textit{910.e1}}, and \href{https://lmfdb.org/EllipticCurve/Q/9438.m2}{\textit{9438.m2}} respectively.
The second and fourth of these are both quadratic twists of the first. It seems that these three have a large Szpiro ratio due to a factor of $3^{21}$ in each of their discriminants. The third has a factor of $2^{63}$ in its discriminant.
These are the only four curves in our set with $ \sigma \ge 8$.
There are $123$ curves in the database with $\sigma \ge 7$, only $15$ of which have conductor larger than~$500000$.
The largest $\sigma$ from our data with $N > 500000$ has $N = 532350$ and $\sigma \approx 7.161459$.

\subsection{Comparison with Cremona's database}
\label{secComparisionToCremonaDB}
Cremona~\cite{Cremona97book,CremonaData} has computed the set of all rational elliptic curves with conductor less than various bounds, currently up to $N\leq 500000$.

If $S$ is the set of primes of bad reduction of an elliptic curve $E$ of conductor $N$, then
\[
N_S\leq N \leq 1728 N_S^2
\]
Thus in principle, the problems of computing $M(S)$ and all curves of bounded conductor are equivalent. Both parameters~$S$ and~$N$ stratify the infinite set of rational elliptic curves.
In practise however these stratifications differ considerably: for example, $M(S(6))$ contains $14216$ curves of conductor $2^8 3^5 5^2 7^2 11^2 13^2 \approx 10^{12}$, which is considerably larger than $500000$; and on the other hand, $M(S(6))$ does not contain the four curves with conductor~$17$.

Cremona's database contains at present $1238682$ distinct $j$-invariants, whereas the computation we performed resulted in $34960$, because for each $j$-invariant, our set contains at least $128$ distinct twists.
On the other hand, Cremona's database contains $3064705$ $\QQ$-isomorphism classes of curves, whereas our contains $4576128$.
Despite the fact that the two databases contain more or less the same number of curves, there are $4376070$ curves in our set not contained in the Cremona database, that is, less than 5\% of our set overlaps with his.
%
We observe significant differences in the distributions of $\log(N)$ and of $\sigma$ for both data sets, see Figures~\ref{figLogConductors} and~\ref{figSzpiroRatios}.

Cremona's tables contain a lot more information about each curve present there than our tables currently do, including Manin constants, generators for the Mordell--Weil group, BSD invariants, modular degrees, optimality data, sets of integral points and image types of Galois representation.
Much of this data would be prohibitively difficult to compute for every curve in our set due in part to the size of the conductors of some of the curves in our table.

\section{Computation}
In this section we discuss the reduction of computing $M(S)$ to the problem of solving Mordell equations, the computation of the requisite Mordell--Weil bases, which is then the dominant computational task to be undertaken, and the heuristic completeness of the obtained data.

The code implementing the methods described here and computed data are available online.
The repository
\href{https://github.com/elliptic-curve-data/ec-data-S6}{\textit{https://github.com/elliptic-curve-data/ec-data-S6}}
contains the majority of the code, and
the file \texttt{mordell.sage} of \href{https://github.com/bmatschke/solving-classical-diophantine-equations/blob/master/mordell.sage}{\textit{https://github.com/bmatschke/solving-classical-diophantine-equations/}} contains an implementation of the algorithm of von K\"anel and the second author~\cite{vKMa14sUnitAndMordellEquationsUsingShimuraTaniyama}.

\subsection{Computation method}
\label{secComputationMethod}
Let $S$ denote a finite set of rational primes, let $M(S)$ denote the set of elliptic curves over $\QQ$ with good reduction outside of~$S$, up to $\QQ$-isomorphism, and let
\[
N_S:=\prod_{p\in S} p\text.
\]

For this section we assume that $2,3\in S$, which can be achieved by enlarging~$S$ if necessary.
Let $\OS=\ZZ[1/N_S]$ denote the ring of $S$-integers and $\OS^*$ the group of $S$-units.

A theorem of Shafarevich~\cite{shafarevich1962algebraic} 
states that for any~$S$ the set of curves $M(S)$ is finite.
This can be seen as follows:
For any $E\in M(S)$ choose a minimal Weierstrass model for $E$ and consider the $c_4$ and $c_6$ invariants and discriminant~$\Delta_E$ of this model.
These invariants satisfy the equation $c_6^2 = c_4^3 - 1728\Delta_E$ and $\Delta_E\in \ZZ\cap \OS^*$.
If necessary we may divide this equation by a power of $p^6$ for each~$p\in S$ to obtain an equality of the form $Y^2 = X^3 + a$, where $X,Y\in \OS$ and $a=\pm \prod_{p\in S} p^{e_p}$ with $0\leq e_p \leq 5$ ($p\in S$).
The pair $(X,Y)$ can then be regarded as an $S$-integral point on the Mordell curve $E_a\colon y^2 = x^3 + a$.
By a theorem of Siegel~\cite{siegel29anwendungenDiophantApprox,Silverman86arithmeticBook}, $E_a(\OS)$ is finite.
From any point in $E_a(\OS)$ we can recover potential invariants $c_4$ and $c_6$ that produce the point, up to any factors of $p^6$ in $c_4^3$ and $c_6^2$ for $p\in S$.
This recovers $E$ up to a quadratic twist by a positive $S$-unit.
Moreover there are exactly $2^{|S|}$ such twists.

We deduce that $M(S)$ is finite and its computation reduces to the computation of~$E_a(\OS)$ for finitely many values of~$a$.
To determine $E_a(\OS)$ we use the algorithm of von K\"anel and the second author~\cite{vKMa14sUnitAndMordellEquationsUsingShimuraTaniyama}, who gave a method to compute $S$-integral points on rational elliptic curves $E$ provided that generators of the free part of~$E(\QQ)$ are known. 
%
%
Their implementation uses an elliptic logarithm sieve, which can compute~$E_a(\OS)$ in quite an efficient manner.
Thus to compute $M(S)$ it turns out that computing the necessary Mordell--Weil bases of $2\cdot 6^{|S|}$ Mordell curves is the computational bottleneck.
In Section~\ref{secMW} we discuss this in detail.

\subsection{Computing Mordell--Weil bases}
\label{secMW}
We have carried out the approach outlined above for $S = S(6)$. We now discuss the most computationally intensive part of the process, which is finding the generators of the free part of the Mordell--Weil group for a number of Mordell curves, many of which have large discriminant.
We will use the term Mordell--Weil basis to refer to these generators.
Note that finding the generators of the torsion subgroup is both computationally easier and completely classified for Mordell curves~\cite{fueter}, so we assume it is known from now on.

The curves we consider are those with
\begin{equation}
\label{eqAs_for_S6}
a \in \{\pm 2^{e_2}3^{e_3}5^{e_5}7^{e_7}11^{e_{11}}13^{e_{13}} \colon 0\le e_p \le 5\},
\end{equation}
giving us $93312$ curves to find the Mordell--Weil bases of.

We can reduce the number of curves that we need to consider using the following fact.

\begin{lemma}\label{lem:threeisog}
    All Mordell curves have a 3-isogeny given by
    \begin{align}
        y^2 = x^3 +a \amp\to y^2 = x^3 - 27a\\
        (x,y) \amp\mapsto \left(\frac{y^2 + 3a}{x^2}, y\frac{y^2 - 9 a}{x^3}\right)
    \end{align}
\end{lemma}

As the composition of two such isogenies is an isomorphism between two models of the same curve, these $3$-isogenies partition our set of Mordell curves into pairs. The upshot is that if we can find generators of the Mordell--Weil group of one of each pair we can easily find generators for the other by pushing the basis forward along the isogeny and saturating, if necessary.
Using this we need only compute the Mordell--Weil bases of half the curves, and we may choose which of each pair to consider.

\subsubsection{Standard techniques}
Out of the $93312$ Mordell curves, we have computed what should be the analytic rank of those with positive $a$ using Pari/GP's \texttt{ellanalyticrank}~\cite{pari2}, via Sage~\cite{sage}.
Using the above isogeny, the other half of the curves will have the same Mordell--Weil rank.
Of these curves, $20215$ have analytic rank~$0$, $23186$ have analytic rank~$1$, $3112$ have analytic rank~$2$, $142$ have analytic rank~$3$, and only $1$ curve has analytic rank~$4$ which is
\[y^2 = x^3 + 82063881900\text.\]

We assume that the output of \texttt{ellanalyticrank} is correct and that the analytic ranks are as stated above. As Pari does not use interval arithmetic it is not clear to what extent these computations are guaranteed to be correct (especially for the high rank cases). As we shall see below we have found as many generators as there should be for almost all curves.
For many of the curves, once a set of generators is found descent techniques can be used to prove that the algebraic rank equals to what is implied by BSD, and that the set of generators is complete.

By the work of Gross--Zagier~\cite{GrossZagier1986heegner} and Kolyvagin~\cite{kolyvagin2007euler} it is known that analytic rank~$\le 1$ implies the rank equals the analytic rank.
Therefore no further computation is required for the analytic rank~$0$ curves above.
For the analytic rank~$1$ curves we need only to find a single non-torsion point which we can then saturate to find a basis.

For many rank~$1$ and~$2$ curves in the set and for all curves of rank at least~$3$, a combination of the built-in Magma and Sage functions and a few other techniques summarized below sufficed to compute the Mordell--Weil bases.
These included two and four-descents methods, point searching with Stoll's \texttt{ratpoints} program~\cite{Stoll14ratpoints} and Simon's \texttt{ellQ}~\cite{Simon02ellQ} to search for points in some instances.

For the curves of rank at least~$2$, sometimes it was only possible to find a subset of a set of generators on each curve of each three-isogenous pair. However in this case it was often possible to mapping one set of generators via the isogeny to the other curve, and combine the generators to give a basis for the Mordell--Weil group of one (and therefore both) curves. This happened mostly when the height of he found generators grew when mapped to the isogenous curve.

In rank~$1$, Heegner points are available in addition to the other machinery of point searching and descent \cite{cohen_number_2007}.
In theory, computing a Heegner point is guaranteed to terminate and if the found point is non-torsion then it is known that the curve has algebraic rank~$1$.
However in order to compute Heegner points we need to find the images of points under the modular parameterization and hence we may need to compute a large number of Frobenius eigenvalues to find the image to a large enough precision in order to recover an algebraic point.

Using a combination of all of these techniques we found bases for all curves but $16481$ of the rank~$1$ curves and we found a single generator (but not the full basis) for all but~$33$ of the rank~$2$ curves. There was one additional rank~$2$ curve for which we did not find any infinite order points with these methods ($E_a$ for $a = 2 \cdot 3 \cdot 5 \cdot 7 \cdot 11^4 \cdot 13^5$). It is likely that a part of the rank~$1$ cases would be amenable to the techniques mentioned, by using larger search bounds or more time or memory. 
However it seemed a different approach was needed to find bases on the hard rank~$1$ curves as well as all remaining rank~$2$ curves.

\subsubsection{$12$-descent}
To determine the generators on these harder curves we used the 12-descent routine in Magma designed and implemented by Fisher~\cite{fisher_finding_2007}.
This works by combining a $3$-cover obtained from a $3$-descent procedure with a $4$-cover from doing $2$-descent and then $4$-descent.
In our setting the presence of a $3$-isogeny for all of our target curves allows us to use $3$-descent by isogeny to obtain the $3$-cover, this is more efficient as the number fields involved are smaller than a general $3$-descent. The implementation for this in Magma is due to Creutz.

Fisher's algorithm then determines a $12$-cover and a map to the original curve from each pair of one $3$-cover and one $4$-cover coming from these lower descents.
Therefore to find a generator of the Mordell curve we loop over all $4$-covers and $3$-covers of the curve coming from descent and search for points on the corresponding $12$-cover.
It is expected that if an $n$-cover has small enough coefficients that the height of a preimage of a point of height $h$ is roughly $h/2n$.
Therefore given an estimate of the canonical height of a generator of the Mordell curve (coming from the regulator estimated via BSD) and a bound for the difference of the na{\"i}ve and canonical heights on an elliptic curve (such as \cite{mullerstoll}) we can search for points on the cover which should be mapping to a generator.
Because this point should have smaller height this should substantially reduce the time needed to search for points, compared with simply searching on the original curve.
Using this we reduce the height to be searched up to by a factor of up to~$24$ if the coefficients of the $12$-cover are not too large.
To search for points on the $12$-covers we use the Magma method \texttt{PointSearch}, implemented by Watkins~\cite{watkinspadic}, see also~\cite{womack_explicit_2003}.
This approach has been used to used to find generators of large height on single Mordell curves previously \cite{weigandt}. 

Due to the fact that we do not know $|\Sha|$ for our curves, the regulator may give an overestimate for the height of a generator, as BSD will only allow us to determine $\sqrt{|\Sha|} \cdot R$ from readily available information.

This procedure was carried out with increasing timeout, up to a maximum of~$12$ hours, and was broadly successful in finding a generator of the rank~$1$ and~$2$ curves for which more standard methods failed.

\subsubsection{Remaining curves}
The combination of these methods has been broadly successful. However there are~$306$ rank~$1$ curves remaining (up to the $3$-isogeny above), which we have so far been unable to find the Mordell--Weil bases of.
A combination of large conductor and large regulator (and hence either large generator height or large $|\Sha|$) has prevented any of the above methods from working in a reasonable time frame.

The Mordell curve with smallest regulator for which we do not know a generator is
\[
y^2 = x^3 + 730033053750
\]
with regulator approximately $167.305352$.

The largest regulators occurring for the remaining curves arise for
\[
    y^2 = x^3 \pm 904509009004500900000,
\]
which interestingly are quadratic twists of each other (by $-1$).
Their regulators are $17550.10$ 
in the $+$ case and $17628.52$ 
in the $-$ case.
However these curves are somewhat exceptional, not all curves are quite so large.
The mean of the remaining regulators is $2622.49$. 

To attack the remaining curves, several options exist to compute a generator.
We have trialled these for a few of the remaining curves that we expect to be ``easier''.

We attempted to make use of Magma's \texttt{HeegnerPoint} method.
As described in Watkins~\cite{watkins_remarks_2006} this allows the user to use $4$-descent to construct a $4$-cover of the target elliptic curve and then find a Heegner point on the cover, reducing the required precision needed and hence the number of required Frobenius traces.
Unfortunately the Magma method fails on many of our difficult examples, presumably because both the conductor of the curve and the height of the Heegner point are both large enough that the number of Frobenius traces needed becomes unwieldy for Magma. Due to the closed source nature of Magma (and the \texttt{HPInternal2} and \texttt{FrobeniusTracesDirect} methods in particular) we have been unable to rectify these problems.
It is also unclear whether or not Magma's algorithm for computing all traces of Frobenius for primes below a given bound for one of our curves is optimal. As our curves are Mordell curves they have CM (by $\sqrt{-3}$), therefore to compute the Frobenius traces we may make use of Cornacchia's algorithm \cite[pp.~597]{cohen_number_2007}.

The highly optimised \texttt{smalljac} package \cite{kedlaya_computing_2008} (available from Sutherland's webpage) includes an implementation of this algorithm in the case of $j$-invariant~$0$ and we expect that using this will be the most effective way to compute enough Frobenius traces to find a Heegner point on the remaining curves.

Happily Pari/GP's \texttt{ellheegner} method is more reliable on our examples, though it does appear to use the covering method, thus we expect that it will conclude on several of the remaining curves given enough time. It is not clear that Cornacchia's algorithm is used to compute all of the Frobenius traces.

For instance this function has returned successfully for one of the ``missing'' curves.
We have found a generator of
\[
    y^2 = x^3 + 4259854045547100000
\]
of height $956.2822$, and it is possible this case was more tractable due to the fact that this is really the double of a generator which has height $239.07055$ instead.
As we are not actually missing a generator on this curve we may check that indeed it does not give any extra elements of $M(S(6))$, however as it took far longer to find this point and it required more interactive experimentation with parameters than the descent methods that we used for the vast majority of the curves we prefer to present it separately to the main data.

In theory with an increased height bound for point searching on 12-covers and with enough time a point should be found on such a cover in the same way as we found the above.
There are two potential issues with this.
Firstly a lattice reduction algorithm is used in the point search procedure. 
It often happens that this method gets stuck if these lattices happen to be ill-conditioned for Magma's algorithm. 
This can stall the point search and we are not aware of the true cause or of ways of avoiding this other than restarting and hoping to get lucky.
The second is that the coefficients of the $12$-cover can be quite large, which can reduce the effectiveness of the height saving of the algorithm.
Thus it is very important to minimise the $12$-cover, as described by Fisher, as much as is possible in order to get the most use out of the method. 
It is plausible that with more work minimizing the $12$-covers the runtime of point searching can be made more feasible.

We have checked another ``missing'' example where $12$-descent succeeds with more individual care than we were able to take at scale. This was curve $E_a$ for $a= 139413405126996000$, which has regulator $1504.24027$, with a height bound of~$10^{21}$ on the $12$-covers the descent finds a point which gives us a generator of height $1504.24027$ on $E_a$.
It is interesting that this point is not a multiple of any smaller generator, suggesting that $\Sha$ is trivial here.

Higher descents are also a potential avenue to complete the process of finding generators for the remaining curves. The work of Fisher allows one to combine covers of coprime degrees subject to a numerical condition on the degrees. This includes the case of combining an $n$ cover with an $n+1$ cover to obtain an $n(n+1)$ cover.
This could conceivably be used to compute $8\cdot 9 = 72$ covers on Mordell curves by combining 8-descent and 9-descent (as a second $p$-isogeny descent) both of which have been implemented in Magma.
It is unknown at present how to make describing and combining such covers practical however.

\subsection{Completeness of the data}
\label{secCompletenessHeuristic}
First, for many Mordell curves we have computed what should be their rank by computing the analytic rank.
This is easier to compute than the algebraic rank in general. According to BSD these ranks are equal, but this is not known in general.
Computing the algebraic rank is more computationally intensive and can be obstructed by non-trivial $\Sha$. However the analogous computation was performed in \cite{vKMa14sUnitAndMordellEquationsUsingShimuraTaniyama} for~$S(5)$.

We have in some cases allowed Magma to assume GRH, which speeds up computation of class groups and hence descent machinery. This does not invalidate searching for points on the corresponding covers, any rational points found are then verified unconditionally to be independent elements of the Mordell--Weil group, but when proving that the algebraic ranks agree with the analytic ones, either GRH or a longer computation time is required.

Secondly and more seriously, we are missing any $S$-integral points on 612 Mordell curves $E_a$ of rank $1$, because so far we were not able to find the generator of Mordell--Weil for 306 curves (as once one curve from each isogeny class's basis is found, the other may be computed relatively easily).
Assuming BSD we may estimate the regulators of these curves up to a factor of $\sqrt{|\Sha|}$.
In the rank~$1$ case the regulator is simply the height of a Mordell--Weil generator. 
So we have that in the missing cases either the generators are of large height or $\sqrt{|\Sha|}$ is large as their product is at least~$150$.

To relate this to the $S$-integral points on these curves, we recall that the \(abc\) conjecture can be used to prove the weak Hall conjecture, which states that integral points $(x,y)$ on the Mordell curve $E_a\colon y^2 = x + a$ satisfy $x = O(a^{2+\eps})$ for any $\eps>0$, see Schmidt~\cite{Schmidt91diophantineApproximationsBook}.
The same proof can be used to show (asymptotic) upper height bounds for $S$-integral points on~$E_a$.
These make it seem unlikely that an $E_a$ of rank~$1$ with a very large Mordell--Weil generator has an $S$-integral point.
These estimates could be made explicit if we assume for example Baker's explicit \(abc\) conjecture~\cite{BakerABC}.
We give more details on this heuristic in Section~\ref{secSHallandABC}.

These missing Mordell--Weil generators of curves of rank $1$ could be computed via the Heegner point method, which is for example implemented in Pari/GP~\cite{pari2}, whose complexity to find $P\in E_a(\QQ)$ is proportional to $\sqrt{N}h(P)$.
Thus together with BSD we estimate that we can prove completeness of our database in about $50$ CPU years.
This is probably less than the (quote) ``many thousand machine hours on 80 cores'' that Bennett, Gherga and Rechnitzer~\cite{BennettGhergaRechnitzer19ellipticCurvesOverQ} used to recompute the database of \cite{vKMa14sUnitAndMordellEquationsUsingShimuraTaniyama} for~$S(5)$.
The original computation of $M(S(5))$~\cite{vKMa14sUnitAndMordellEquationsUsingShimuraTaniyama} was not timed, but recalling from memory it took in the order of one CPU year.

\subsection{An $S$-integral weak Hall conjecture and the $abc$ conjecture}
\label{secSHallandABC}
In this section we will discuss an $S$-integral analogue of the classical Hall conjecture and how it adds to our heuristic for why our database should be complete.
As for the classical Hall conjecture, we will show that it is implied by the $abc$ conjecture.

For this section we will use the following terminology. 
For any finite set of rational primes~$S$, we call a pair of integers $(x,y)$ \textdef{$S$-primitive} if there is no $p\in S$ such that $p^6$ divides both $x^3$ and~$y^2$.
We formulate an $S$-integral generalization of the weak Hall conjecture.

\begin{conjecture}[An $S$-integral weak Hall conjecture]
\label{conjSHall}
Let $S$ be a finite set of rational primes.
Let $D\neq 0$ be an integer.
For any $\eps>0$, any $S$-primitive solution $(x,y)$ of the equation
\begin{equation}
\label{eqSintegralHallEquation}
y^2 = x^3 + aD,\qquad x,y\in \ZZ, \quad a\in\ZZ\cap\OS^\times,
\end{equation}
satisfies
\begin{equation}
\label{eqSintegralWeakHallInequality}
\max(|x|^{1/2},|y|^{1/3}) = O_\eps((N_SD)^{1+\eps}).
\end{equation}
\end{conjecture}

Recall that the $abc$ conjecture states that for any $\eps>0$ the following holds.
If $a,b,c$ are coprime integers with $a+b+c=0$, then 
\begin{equation}
\label{eqABC}
\max(|a|,|b|,|c|)\leq O_\eps(\rad(abc)^{1+\eps}),
\end{equation}
where $\rad(abc) = \prod_{p\divides abc} p$.
More explicitly, \eqref{eqABC} states that $\max(|a|,|b|,|c|)\leq K_\eps\,\rad(abc)^{1+\eps}$, where $K_\eps$ is a constant that depends only on~$\eps$.

\begin{theorem}
The $abc$ conjecture implies the $S$-integral weak Hall conjecture.
More explicitly, if the $abc$ conjecture holds for some $0<\eps\leq 0.1$ with constant $K_\eps$, then any $S$-primitive solution $(x,y)$ of~\eqref{eqSintegralHallEquation} satisfies
\begin{equation}
\label{eqWeakSintegralHallBound_from_ABC}
\max(|x|^{1/2},|y|^{1/3}) \leq K_\eps^{1+10\eps}(N_sD)^{1+12\eps}.
\end{equation}
\end{theorem}

Our proof largely follows  Schmidt's proof~\cite{Schmidt91diophantineApproximationsBook} that $abc$ implies the classical weak Hall conjecture, although the proof below avoids some technicalities by choosing $s$ and $t$ (see proof) in an efficient way.

\begin{proof}
Suppose $(x,y)$ is an $S$-primitive solution of~\eqref{eqSintegralHallEquation}.
Let $g=\gcd(x^3,y^2)$.
Let $A = x^3/g$, $B = -y^2/g$ and $C=aD/g$, which are coprime integers.
As $A+B+C=0$, the $abc$ conjecture implies that
\begin{equation}
\label{eqabcForABCinProof}
\max(|x|^3/g,|y|^2/g) \leq K_\eps\, \rad(ABC)^{1+\eps}.
\end{equation}
We claim that 
\begin{equation}
\label{eqRadABCdividesXYNSDoverG}
\rad(ABC) \divides \frac{xyN_SD}{g}.
\end{equation}
To see this we consider two cases.

Case 1.) If some $p\in S$ divides $ABC$, then by $S$-primitivity of $(x,y)$ we have $\ord_p(y)\leq 2$ or $\ord_p(x)\leq 1$.
In either case, $\ord_p(g)\leq 4$.
If $\ord_p(g)=4$, then $p^2\divides x$, $p^2\divides y$, $p\divides N_S$, and thus $p\divides xyN_SD/g$.
The cases $\ord_p(g)\in\{0,2,3\}$ are similar, and $\ord_p(g)=1$ is a priori not possible.

Case 2.) Suppose some $p\not\in S$ divides $ABC$.
If $p\nmid g$ then the obvious $p\divides xyD$ suffices.
If $p\divides g$, then $\ord_p(xyD/g)\geq \ord_p(g)(1/3+1/2+1-1) >0$ and so $p\divides xyD/g$.
This finishes the proof of~\eqref{eqRadABCdividesXYNSDoverG}.

Plugging \eqref{eqRadABCdividesXYNSDoverG} into~\eqref{eqabcForABCinProof} implies that
$\max(|x|^3,|y|^2)\leq K_\eps (xyN_SD)^{1+\eps}$ and hence
\[
|x|^{3s}|y|^{2t}\leq K_\eps^{s+t}(xyN_SD)^{(s+t)(1+\eps)}.
\]
For $s = (1-\eps)/(1-5\eps)$ and $t = (1+\eps)/(1-5\eps)$ we obtain
\[
|x| \leq K_\eps^{2/(1-5\eps)} (N_SD)^{(2+2\eps)/(1-5\eps)}.
\]
Similarly for $s = (1+\eps)/(1-5\eps)$ and $t = (2-\eps)/(1-5\eps)$ we obtain
\[
|y| \leq K_\eps^{3/(1-5\eps)} (N_SD)^{(3+3\eps)/(1-5\eps)}.
\]
This yields
\[
\max(|x|^{1/2},|y|^{1/3}) \leq K_\eps^{1/(1-5\eps)} (N_SD)^{(1+\eps)/(1-5\eps)}.
\]
For $\eps\leq 0.1$ this reduces to the claimed bounds.
\end{proof}

Let us relate this to $S$-integral points on the above Mordell curves~$E_a\colon y^2 = x^3 + a$, where $a$ is as in~\eqref{eqAs_for_S6} an $S$-unit with bounded exponents.
Suppose $P=(X,Y)\in E_a(\OS)$.
We can clear denominators of $X$ and $Y$ by multiplying $X^3$ and $Y^2$ by suitable powers of $p^6$ for each $p\in S$, and call the resulting integers $\wt X$ and $\wt Y$.
This yields a relation $\wt Y^2 = \wt X^3 + \wt a$, to which we can apply the $S$-integral weak Hall conjecture~\ref{conjSHall} (with $D=1$), or alternatively~\eqref{eqWeakSintegralHallBound_from_ABC} as implied by the $abc$ conjecture.
We obtain conjectural asymptotic height bounds for $|\wt X|^3$ and $|\wt Y|^2$, which imply up to a small explicit constant (depending on~$S$) the same bound on the na\"i{}ve height of $P$, which in turn is up to an explicitly bounded error the N\'eron--Tate height~$\hat h(P)$. 

In case $S=S(6)$ we can thus make the following heuristic.
First, assume that the $abc$ conjecture holds for $\eps=0.1$ with a constant $K_\eps\leq 1.1\cdot 10^8$.
We checked that this bound indeed holds for all $abc$-triples of the ABC@Home project by de~Smit~\cite{deSmitABCatHome} for which we could compute the radical. 
Using this $\eps$ and $K_\eps$ and the above reasoning, we would obtain a bound for $\hat h(P)$ of approximately $2(2\log K_\eps + 2.2\log N_S)\leq 120$.


\newcommand{\Nmax}{M_S}

\section{Attainability of maximal conductor by curves in $M(S)$}
\label{secMaximalConductor}
In this section, we prove some results suggested by empirical observations of our data.

Specifically we ask the following: Given a set of rational primes~$S$, what is the highest possible conductor of an elliptic curve over $\QQ$ with good reduction outside~$S$?
An immediate upper bound is $\Nmax:= \prod_{p\in S} p^{f_p}$ where $f_2=8$, $f_3=5$, and $f_p = 2$ for $p\geq 5$.
More specifically we may then ask:

\begin{question}
\label{quMaximalConductor}
Does there exist a curve of conductor of $\Nmax$ for any set $S$?
\end{question}

The answer to this question is know without further conditions on~$S$.
For example there does not exist an elliptic curve with good reduction away from~$5$, however the answer is positive for a large class of $S$.
Motivated by our data we have the following sufficient criterion.

\begin{theorem}
\label{thmMaximalConductor}
Let $S$ be a finite set of rational primes that contains either $2$ or $3$ (or both).
Then there exists an elliptic curve over $\QQ$ with conductor $N=\Nmax$.
\end{theorem}

In order to prove the theorem we recall the notion of quadratic twists of elliptic curves. 
For any rational elliptic curve $E\colon  y^2 = x^3 + ax + b$ and an integer $d$, we denote by $E^d\colon  y^2 = x^3 + d^2ax + d^3b$ its quadratic twist by~$d$. 

The theorem now follows immediately from the following lemma.
The proof is constructive.

\begin{lemma}
\label{lemMaximalConductor}
Let $d$ be a square-free product of primes $p\geq 5$.
\begin{enumerate}
\item
Let $E_{\{2,3\}}\colon y^2 = x^3 - 18x + 24$.
Then $E_{\{2,3\}}^d$ has conductor $N = 2^8  3^5 d^2$ and Kodaira type $\mathrm{III}$ at $2$, $\mathrm{II}$ at $3$, and $\mathrm{I}_0^*$ at $p\geq 5$ with $p\divides d$.
\item
Let $E_{\{2\}}\colon y^2 = x^3 + 8x$.
Then $E_{\{2\}}^d$ has conductor $N = 2^8 d^2$ and Kodaira type $\mathrm{III}^*$ at $2$ and $\mathrm{I}_0^*$ at $p\geq 5$ with $p\divides d$.
\item
Let $E_{\{3\}}\colon  y^2 + y = x^3 - 1$.
Then $E_{\{3\}}^d$ has conductor $N = 3^5 d^2$ and Kodaira type $\mathrm{II}$ at $3$ and $\mathrm{I}_0^*$ at $p\geq 5$ with $p\divides d$.
\end{enumerate}
\end{lemma}

\begin{proof}
This is a straightforward computation with Tate's algorithm, which we omit here.
For the convenience of the reader it is available as an appendix of the GitHub and arXiv version of this paper, which can be found at:

\href{https://github.com/elliptic-curve-data/ec-data-S6/blob/master/docs/paper.pdf}{\textit{https://github.com/elliptic-curve-data/ec-data-S6/blob/master/docs/paper.pdf}}
\end{proof}

We remark that in general, twisting an elliptic curve~$E$ by a prime $p\geq 5$ may change the reduction type of~$E$ at~$2$ and~$3$, but this does not happen for the three curves listed in the lemma.

Silverman~\cite[Exercises 4.52, 4.53]{Silverman94advancedTopicsBook} gives two families of elliptic curves defined over~$\QQ$, which have maximal possible conductor exponent at~$3$ and at~$2$, respectively, and also have this property after base changing to a number field.
The above curve $E_{\{2\}}$ belongs to Silverman's latter family.

\section{Applications}
\label{secApplications}

In this section we will briefly discuss some applications of the dataset.

\subsection{Solving $S$-unit equations}
\label{secSUE_connection}

Let $S$ be a finite set of rational primes. 
As above denote by $\OS$ and $\OS^*$ the $S$-integers and $S$-units, respectively.
The $S$-unit equation is the equation
\begin{equation}
\label{eqSunit}
x+y=1, \qquad x,y\in\OS^*.
\end{equation}
This classical diophantine equations is intimately related to the \(abc\) conjecture, this can be seen by clearing denominators to obtain an $abc$ equation.
Also, more generally, $S$-unit equations over number fields are known to have only finitely many solutions, as was first shown by Siegel~\cite{siegel29anwendungenDiophantApprox} 
and Mahler~\cite{mahler1933approximation}. 
Siegel~\cite{siegel29anwendungenDiophantApprox,Silverman86arithmeticBook} used this to prove that any hyperelliptic curve of genus at least one has only finitely many $S$-integral points.

It turns out that solving $S$-unit equations can be reduced to the computation of $M(S\cup\{2\})$ via Frey--Hellegouarch curves:
If $(x,y)$ is a solution of the $S$-unit equation, then $E_x\colon Y^2 = X(X-1)(X-x)$ lies in $M(S\cup\{2\})$.
Moreover any curve $E\in M(S\cup\{2\})$ can be obtained in this way from at most six different solutions of~\eqref{eqSunit}, and these can be computed from the six possible modular $\lambda$-invariants of~$E$.
In our case, \eqref{eqSunit} for $S=S(6)$ is exactly the case that has been considered by de Weger~\cite{deWeger87solvingExponentialDiophantineEquations}.
He proved that, up to symmetry, it has exactly 545 solutions.
We checked that the curves associated to all of these can be found in our database, which means that our database certainly contains all Frey--Hellegouarch curves with good reduction outside~$S(6)$.
We remark that \eqref{eqSunit} has been solved for~$S=S(16)$, as well as for all $S$ with $N_S\leq 10^7$~\cite{vKMa14sUnitAndMordellEquationsUsingShimuraTaniyama}. This is far out of reach for the above method of reducing~\eqref{eqSunit} to computing~$M(S)$.

In the other direction, the computation of~$M(S)$ can be reduced to solving $S'$-unit equations over finitely many number fields, where the number fields are all possible number fields $K$ of degree at most six that are unramified outside $S\cup\{2\}$ and $S'$ being the primes in~$K$ above $S\cup\{2\}$. This link was made into an algorithm by Koutsianas~\cite{Koutsianas19ellipticCurvesOverNFs}.

\subsection{Other diophantine problems}
\label{secReductionsOfDiphantineEquations}
Many other diophantine problems reduce to the computation of~$M(S)$, notably cubic Thue--Mahler equations
\[
ax^3 + bx^2y + cxy^2 + dy^3 = m\prod_{p\in S}p^{e_p}, \qquad x,y\in\ZZ,\ \ e_p\in\ZZ_{\geq 0} \ (p\in S),
\]
where $a,b,c,d,m \in\ZZ$ and $m\neq 0$ are given such that the left-hand side has non-vanishing discriminant. Likewise generalized Ramanujan--Nagell equations
\[
x^2+b = y, \qquad x\in\OS,\  y\in\OS^*,
\]
where $b\neq 0$ is a given integer, can be reduced to computing~$M(S)$.
In particular we can find solutions for these equations for $S=S(6)$ via our computation of curves in $M(S)$, which subject to the hypothesis that we have in fact found the whole set $M(S)$ should be the complete sets of solutions of these equations, see the above discussion on completeness in Section~\ref{secCompletenessHeuristic}.

\subsection{$n$-congruences between elliptic curves}
Given $n\in \NN$, a pair of elliptic curves $E_1,E_2/\QQ$ for which $E_1[n]\simeq E_2[n]$ as Galois modules are called $n$-congruent. The Frey--Mazur conjecture implies that there should be an absolute bound $C$ such that if $p\ge C$ and $E_1,E_2$ are $p$-congruent, then $E_1$ and $E_2$ must be isogenous.
The only known example of a pair of non-isogenous 17-congruent elliptic curves, found by Cremona and then Billerey \cite{billerey2016remarkable}, occurs for a pair of curves with good reduction outside of ${3, 5, 7, 13}$.
Using our database we searched for similar examples of $n$-congruences between curves for primes $13\le n\le 47$. We found several instances of 13-congruences that were outside the range of existing databases. Fisher has recently found an infinite family of 13-congruent curves \cite{fisher13}, which the examples in our database are all members of.
We did not find any further examples of 17 (or higher) congruences between curves in our database,  other than quadratic twists of the example of Cremona--Billerey mentioned above.

\appendix

\section{Proof of Lemma~\ref{lemMaximalConductor}}

1). 
Let $E:=E_{\{2,3\}}:  y^2 = x^3 - 18 x + 24$, which is a minimal model with $j = 5184 = 2^6 3^4$ and (minimal) discriminant
$\Delta(E) = 2^9 3^5$.
$E$ has Kodaira type~$\mathrm{III}$ at~$2$, and $\mathrm{II}$ at~$3$.
Its conductor is thus~$2^8 3^5$.
Now let~$d$ be a (square-free) product of primes that are at least~$5$.
Let $E^d: y^2 = x^3 - 18 d^2 x + 24 d^3$ be the quadratic twist of $E$ by~$d$.
$E^d$ has discriminant $\Delta(E^d)=2^9 3^5 d^6$, 
 $a_2 = 0$,
 $a_3 = 0$, 
 $a_4 = -18 d^2$,
 $a_6 = 24 d^3$,
 $b_2 = 0$,
 $b_4 = -36 d^2$,
 $b_6 = 96 d^3$, and
 $b_8 = -324 d^4$.
 First we consider $p \mid d$.
  Let $\wt{d}:=d/p$.
  In Tate's algorithm, we arrive in Step~$6$, and have
  $a_{2,1} = 0$,
  $a_{4,2} = -18 \wt{d}^2$,
  $a_{6,3} = 24 \wt{d}^3$, and 
  $P(t) = t^3 -18 \wt{d}^2 t + 24 \wt{d}^3$.
  $P$ has discriminant $2^5 3^5 \wt{d}^6$, which is non-zero mod~$p$.
  Thus the Kodaira type of $E^d$ at $p \mid d$ is $\mathrm{I_0^*}$.
  Hence $f_p = 2$ (additive reduction at~$p$). 
  
 Next we consider Tate's algorithm for $E^d$ at $p=2$:
  We arrive at Step~$4$ and see that $2^3$ does not divide~$b_8$.
  Thus $E^d$ has Kodaira type~$\mathrm{III}$ at~$2$ and $f_2 = \ord_2(\Delta(E^d))-1 = 8$.
 
 Next we consider Tate's algorithm for $E^d$ at $p=3$:
  We arrive at Step~$3$ and see that $3^2$ does not divide~$a_6$.
  Thus $E^d$ has Kodaira type~$\mathrm{II}$ at~$3$ and $f_3 = \ord_3(\Delta(E^d)) = 5$.

\medskip

2.)
Let $E_2:=E_{\{2\}}:  y^2 = x^3 + 8 x$, which is a minimal model with $j = 1728$ and (minimal) discriminant~$-2^{15}$.
$E_2$ has Kodaira type~$\mathrm{III^*}$ at~$3$.
Its conductor is thus~$2^8$.
Now let $d$ be a (square-free) product of primes that are at least~$5$.
Let $E_2^d: y^2 = x^3 + 8 d^2 x$ be the quadratic twist of $E_2$ by~$d$.
$E_2^d$ has discriminant $\Delta(E_2^d)=-2^{15} d^6$,  
 $a_1 = 0$,
 $a_2 = 0$,
 $a_3 = 0$, 
 $a_4 = 8 d^2$,
 $a_6 = 0$,
 $b_2 = 0$,
 $b_4 = 16 d^2$,
 $b_6 = 0$, and
 $b_8 = -64 d^4$.
 First we consider $p \mid d$.
  Let $\wt{d}:=d/p$.
  In Tate's algorithm, we arrive in Step~$6$, and have
  $a_{2,1} = 0$,
  $a_{4,2} = 8 \wt{d}^2$,
  $a_{6,3} = 0$, and $P(t) = t^3 + 8 \wt{d}^2 t$.
  $P$ has discriminant $-2^{11} \wt{d}^6$, which is non-zero mod~$p$.
  Thus the Kodaira type of $E_2^d$ at $p \mid d$ is~$\mathrm{I_0^*}$.
  Hence for $E_2^d$, $f_p = 2$ (additive reduction at~$p$). 
  
 Next we consider Tate's algorithm for $E_2^d$ at $p=2$:
  We arrive in Step~$6$, and have
  $P(t) = t^3 + 2 d^2 + 2 d^3 = t^3$ mod~$2$.
  P has triple root at $t=0$.
  We arrive at Step~$8$ and consider
  $a_{3,2} = 0$ and $a_{6,4} = 0$.
  So $Y^2 - a_{3,2} Y + a_{6,4}$ has double root at~$0$.
  We arrive at Step~$9$, no coordinate change needed, but $2^4$ does not divide~$a_4$.
  Thus $E_2^d$ has type~$\mathrm{III^*}$ at~$2$ and $f_2 = \ord_2(\Delta(E_2^d))-7 = 8$.)
  
 Next we consider Tate's algorithm for $E_2^d$ at $p=3$:
  We have $\ord_3(\Delta(E_2^d))=0$.
  Thus $E_2^d$ has Kodaira type~$\mathrm{I_0}$ at~$3$ and $f_3 = 0$.

\medskip

3.) 
Let $E_3:=E_{\{3\}}:  y^2 + y = x^3 - 1$, which is a minimal model with $j = 0$, and (minimal) discriminant $-3^5$.
It is isomorphic to $E'_3: y^2 = x^3 - 48$, which has discriminant $-2^{12} 3^5$.
$E_3$ has Kodaira type~$\mathrm{II}$ at~$3$.
Its conductor is thus~$3^5$.
Now let $d$ be a (square-free) product of primes that are at least~$5$.
Let ${E'_3}^d: y^2 = x^3 - 24 d^3$ be the quadratic twist of $E'_3$ by~$d$.
${E'_3}^d$ has
 discriminant $\Delta({E'_3}^d)=-2^{12} 3^5 d^6$ ($p$-minimal for~$p\geq 3$), 
 $a_2 = 0$,
 $a_3 = 0$, 
 $a_4 = 0$,
 $a_6 = -2^4 3 d^3$,
 $b_2 = 0$,
 $b_4 = 0$,
 $b_6 = -2^6 3 d^3$, and
 $b_8 = 0$.
 First we consider $p \mid d$.
  Let $\wt{d}:=d/p$.
  In Tate's algorithm, we arrive in Step~$6$, and have
  $a_{2,1} = 0$,
  $a_{4,2} = 0$,
  $a_{6,3} = -2^4 3 \wt{d}^3$, and 
  $P(t) = t^3 - 48 \wt{d}^3$.
  $P$ has discriminant $-2^8 3^5 \wt{d}^6$, which is non-zero mod~$p$.
  Thus the Kodaira type of ${E'_3}^d$ (and hence of~$E_3^d$) at $p \mid d$ is~$\mathrm{I_0^*}$.
  Hence for $E_3^d$, $f_p = 2$ (additive reduction at~$p$). 

 Next we consider Tate's algorithm for $E_3^d$ at~$p=2$:
  Since $\Delta(E_3^d) = -3^5 d^2$, $\ord_2(\Delta(E_3^d)) = 0$, 
  hence Kodaira type~$\mathrm{I_0}$ at~$2$ and $f_2 = 0$.

 Next we consider Tate's algorithm for ${E'_3}^d$ at $p=3$:
  We arrive at Step~$3$ and see that $3^2$ does not divide~$a_6$.
  Thus $E_3^d$ has Kodaira type~$\mathrm{II}$ at~$3$ and $f_3 = \ord_3(\Delta({E'_3}^d)) = 5$.
%
This finishes the proof of Lemma~\ref{lemMaximalConductor}
\qed

\small

\bibliographystyle{plain}

\bibliography{bib}

\begin{thebibliography}{10}

\bibitem{AgrawalCoatsHuntVDPoorten80}
M.~K. Agrawal, J.~H. Coates, D.~C. Hunt, and A.~J. van~der Poorten.
\newblock Elliptic curves of conductor {$11$}.
\newblock {\em Math. Comp.}, 35(151):991--1002, 1980.

\bibitem{BakerABC}
Alan Baker.
\newblock Experiments on the abc-conjecture.
\newblock {\em Publ. Math. Debrecen}, 65(3-4):253--260, 2004.

\bibitem{BennettGhergaRechnitzer19ellipticCurvesOverQ}
Michael~A. Bennett, Adela Gherga, and Andrew Rechnitzer.
\newblock Computing elliptic curves over {$\QQ$}.
\newblock {\em Math. Comp.}, 88(317):1341--1390, 2019.

\bibitem{BennettRechnitzer}
Michael~A. Bennett and Andrew Rechnitzer.
\newblock Computing elliptic curves over {$\mathbb Q$}: bad reduction at one
  prime.
\newblock In {\em Recent progress and modern challenges in applied mathematics,
  modeling and computational science}, volume~79 of {\em Fields Inst. Commun.},
  pages 387--415. Springer, New York, 2017.

\bibitem{billerey2016remarkable}
Nicolas Billerey.
\newblock On some remarkable congruences between two elliptic curves, 2016.
\newblock arXiv:1605.09205.

\bibitem{BirchKuyk}
Bryan~J. Birch and Willem Kuyk, editors.
\newblock {\em Modular functions of one variable. {IV}}.
\newblock Lecture Notes in Mathematics, Vol. 476. Springer-Verlag, Berlin-New
  York, 1975.

\bibitem{BCDT01modularityOverQ}
Christophe Breuil, Brian Conrad, Fred Diamond, and Richard Taylor.
\newblock On the modularity of elliptic curves over {$\mathbf Q$}: wild 3-adic
  exercises.
\newblock {\em J. Amer. Math. Soc.}, 14(4):843--939, 2001.

\bibitem{brumermcguinness}
Armand Brumer and Oisín McGuinness.
\newblock The behavior of the {Mordell}-{Weil} group of elliptic curves.
\newblock {\em Bull. Am. Math. Soc.}, 23(2):375--382, 1990.

\bibitem{Coghlan67ellipticCurves23}
Francis Coghlan.
\newblock Elliptic {C}urves with {C}onductor $2^m 3^n$.
\newblock Ph.D. thesis, Manchester, England, 1967.

\bibitem{cohen_number_2007}
Henri Cohen.
\newblock {\em Number {Theory}: {Volume} {I}: {Tools} and {Diophantine}
  {Equations}}.
\newblock Graduate {Texts} in {Mathematics}, {Number} {Theory}.
  Springer-Verlag, New York, 2007.

\bibitem{CremonaData}
John~E. Cremona.
\newblock Elliptic curve data.
\newblock
  {\href{https://johncremona.github.io/ecdata/}{https://johncremona.github.io/ecdata/}}.

\bibitem{Cremona97book}
John~E. Cremona.
\newblock {\em Algorithms for modular elliptic curves}.
\newblock Cambridge University Press, second edition, 1997.

\bibitem{CremonaLingham07ellipticCurves}
John~E. Cremona and Mark~P. Lingham.
\newblock Finding all elliptic curves with good reduction outside a given set
  of primes.
\newblock {\em Experiment. Math.}, 16(3):303--312, 2007.

\bibitem{deSmitABCatHome}
Bart de~Smit.
\newblock {ABC@Home}.
\newblock
  {\href{https://www.math.leidenuniv.nl/~desmit/abc/}{https://www.math.leidenuniv.nl/$\sim$desmit/abc/}}.

\bibitem{deWeger87solvingExponentialDiophantineEquations}
Benjamin M.~M. de~Weger.
\newblock Solving exponential {D}iophantine equations using lattice basis
  reduction algorithms.
\newblock {\em J. Number Theory}, 26(3):325--367, 1987.

\bibitem{sage}
The~Sage Developers.
\newblock {\em {S}age{M}ath ({V}ersion 9.0)}, 2020.
\newblock {\href{http://www.sagemath.org}{http://www.sagemath.org}}.

\bibitem{fisher_finding_2007}
Tom Fisher.
\newblock Finding rational points on elliptic curves using 6-descent and
  12-descent.
\newblock \href{https://arxiv.org/abs/0711.3774}{arXiv:0711.3774}, 2007.

\bibitem{fisher13}
Tom Fisher.
\newblock On families of 13-congruent elliptic curves, 2019.
\newblock arXiv:1912.10777.

\bibitem{fueter}
Rudolf Fueter.
\newblock Ueber kubische diophantische {G}leichungen.
\newblock {\em Comment. Math. Helv.}, 2(1):69--89, 1930.

\bibitem{GrossZagier1986heegner}
Benedict~H. Gross and Don~B. Zagier.
\newblock Heegner points and derivatives of $l$-series.
\newblock {\em Invent. Math.}, 84(2):225--320, 1986.

\bibitem{weigandt}
Jamie~Weigandt (https://mathoverflow.net/users/4872/jamie weigandt).
\newblock Rational points on $y^2=x^3-86069^5$.
\newblock MathOverflow.
\newblock https://mathoverflow.net/q/105591 (version: 2012-09-25).

\bibitem{vKMa14sUnitAndMordellEquationsUsingShimuraTaniyama}
Rafael {\noopsort{K\"anel}{von K\"anel}} and Benjamin Matschke.
\newblock {S}olving {$S$}-unit, {M}ordell, {T}hue, {T}hue--{M}ahler and
  generalized {R}amanujan--{N}agell equations via {S}himura--{T}aniyama
  conjecture.
\newblock {\href{http://arxiv.org/abs/1605.06079}{arXiv:1605.06079}},
  {\href{https://bmatschke.github.io/solving-classical-diophantine-equations/}{https://bmatschke.github.io/solving-classical-diophantine-equations/}},
  2016.

\bibitem{kedlaya_computing_2008}
Kiran~S. Kedlaya and Andrew~V. Sutherland.
\newblock Computing {L}-{Series} of {Hyperelliptic} {Curves}.
\newblock In Alfred~J. van~der Poorten and Andreas Stein, editors, {\em
  Algorithmic {Number} {Theory}}, Lecture {Notes} in {Computer} {Science},
  pages 312--326, Berlin, Heidelberg, 2008. Springer.

\bibitem{kolyvagin2007euler}
Victor~A. Kolyvagin.
\newblock Euler systems.
\newblock In {\em The Grothendieck Festschrift}, pages 435--483. Springer,
  2007.

\bibitem{Koutsianas19ellipticCurvesOverNFs}
Angelos Koutsianas.
\newblock Computing all elliptic curves over an arbitrary number field with
  prescribed primes of bad reduction.
\newblock {\em Exp. Math.}, 28(1):1--15, 2019.

\bibitem{lmfdb}
The {LMFDB Collaboration}.
\newblock The {L}-functions and modular forms database.
\newblock \href{https://www.lmfdb.org}{https://www.lmfdb.org}, 2020.

\bibitem{mahler1933approximation}
Kurt Mahler.
\newblock {Zur Approximation algebraischer Zahlen. I}.
\newblock {\em Math. Ann.}, 107(1):691--730, 1933.

\bibitem{mullerstoll}
J.~Steffen Müller and Michael Stoll.
\newblock Computing {Canonical} {Heights} on {Elliptic} {Curves} in
  {Quasi}-{Linear} {Time}.
\newblock {\em LMS J. Comput. Math.}, 19(A):391--405, 2016.
\newblock arXiv:1509.08748.

\bibitem{Ogg66_2powerConductor}
Andrew~P. Ogg.
\newblock Abelian curves of 2-power conductor.
\newblock {\em Math. Proc. Camb. Philos. Soc.}, 62(2):143--148, 1966.

\bibitem{Schmidt91diophantineApproximationsBook}
Wolfgang~M. Schmidt.
\newblock {\em Diophantine approximations and {D}iophantine equations}, volume
  1467 of {\em Lecture Notes in Mathematics}.
\newblock Springer-Verlag, Berlin, 1991.

\bibitem{shafarevich1962algebraic}
Igor~R. Shafarevich.
\newblock Algebraic number fields.
\newblock In {\em Proceedings of an International Congress on Mathematics,
  Stockholm}, pages 163--176, 1962.

\bibitem{siegel29anwendungenDiophantApprox}
Carl~L. {Siegel}.
\newblock {\"Uber einige Anwendungen diophantischer Approximationen.}
\newblock {\em {Abh. Preu{\ss}. Akad. Wiss., Phys.-Math. Kl.}}, 1929(1):70,
  1929.

\bibitem{Silverman86arithmeticBook}
Joseph~H. Silverman.
\newblock {\em Arithmetic of elliptic curves.}, volume 106 of {\em Graduate
  Texts in Mathematics}.
\newblock Springer, 1986.

\bibitem{Silverman94advancedTopicsBook}
Joseph~H. Silverman.
\newblock {\em Advanced topics in the arithmetic of elliptic curves}, volume
  151 of {\em Graduate Texts in Mathematics}.
\newblock Springer, 1994.

\bibitem{Simon02ellQ}
Denis Simon.
\newblock Computing the rank of elliptic curves over number fields.
\newblock {\em LMS J. Comput. Math.}, 5:7--17, 2002.
\newblock Program ellQ:
  \href{https://simond.users.lmno.cnrs.fr/ellQ.gp}{https://simond.users.lmno.cnrs.fr/ellQ.gp}.

\bibitem{SteinWatkins}
William~A. Stein and Mark Watkins.
\newblock A {Database} of {Elliptic} {Curves} — {First} {Report}.
\newblock In Claus Fieker and David~R. Kohel, editors, {\em Algorithmic
  {Number} {Theory}}, Lecture {Notes} in {Computer} {Science}, pages 267--275,
  Berlin, Heidelberg, 2002. Springer.

\bibitem{Stephens65thesis}
Nelson~M. Stephens.
\newblock The {B}irch {S}winnerton-{D}yer {C}onjecture for {S}elmer curves of
  positive rank.
\newblock Ph.D. Thesis, Manchester, 1965.

\bibitem{Stoll14ratpoints}
Michael Stoll.
\newblock Documentation for the ratpoints program.
\newblock \href{https://arxiv.org/abs/0803.3165}{arXiv:math/0803.3165}, 2014.

\bibitem{pari2}
{The PARI~Group}, Univ. Bordeaux.
\newblock {\em {PARI/GP version 2.11.2}}, 2019.
\newblock
  \href{http://pari.math.u-bordeaux.fr/}{http://pari.math.u-bordeaux.fr/}.

\bibitem{Tingley75thesis}
Dave~J. Tingley.
\newblock Elliptic curves uniformized by modular functions.
\newblock Ph.D. thesis, University of Oxford, 1975.

\bibitem{watkinspadic}
Mark Watkins.
\newblock Searching for points using the {E}lkies {ANTS-IV} algorithm.
\newblock
  {\href{http://magma.maths.usyd.edu.au/~watkins/papers/padic.ps}{http://magma.maths.usyd.edu.au/~watkins/papers/padic.ps}}.

\bibitem{watkins_remarks_2006}
Mark Watkins.
\newblock Some {Remarks} on {Heegner} {Point} {Computations}.
\newblock \href{https://arxiv.org/abs/math/0506325}{arXiv:math/0506325}, 2006.

\bibitem{womack_explicit_2003}
Thomas Womack.
\newblock {\em Explicit {Descent} on {Elliptic} {Curves}}.
\newblock PhD thesis, University of Nottingham, July 2003.

\end{thebibliography}

\end{document}